\documentclass[10pt,preprint,reqno]{amsart}
\usepackage{amsmath,amssymb,amsfonts,amsthm}

\usepackage{amscd, calc, tikz, tikz-cd}
\usepackage[all]{xy}
\usepackage[color=pink,textsize=footnotesize]{todonotes}

\usepackage{enumerate,graphicx,psfrag}
\usepackage{mathrsfs}
\usepackage{hyperref}

\usepackage{color} %include even if images aren’t in color
\usepackage{pstricks,pdftricks} 

\usepackage[shortlabels]{enumitem}

\usepackage{cases}  
\usepackage{tikz-feynman}
\usetikzlibrary{decorations.pathmorphing}         

\usepackage[english]{babel}

\usepackage{caption} 
\usepackage{bm}

\newtheorem{thm}{Theorem}[section]
\newtheorem{lem}[thm]{Lemma}

\newtheorem{cor}[thm]{Corollary}

\theoremstyle{definition}
\newtheorem{defn}[thm]{Definition}
\newtheorem{nota}[thm]{Notation}

\newtheorem{example}[thm]{Example}

\allowdisplaybreaks

\newcommand\ds\displaystyle
\newcommand\ts\textstyle

\renewcommand{\phi}{\varphi}                 % Personal preferences.
\renewcommand{\epsilon}{\varepsilon}
\newcommand\eset{\varnothing}

\renewcommand\emptyset\eset
\renewcommand\ell{l}

\newcommand\cupdot {\mbox{\hspace{.15em}$\cup$\hspace{-.47em}$\cdot$\hspace{.4em}}}

\newcommand\supp{\operatorname{supp}}

\newcommand\cB{\mathcal{B}}
\newcommand\cC{\mathcal{C}}
\newcommand\cD{\mathcal{D}}

\newcommand\cM{\mathcal{M}}
\newcommand\cN{\mathcal{N}}

\newcommand\cU{\mathcal{U}}
\newcommand\cV{\mathcal{V}}

\newtheoremstyle{case}{}{}{}{}{}{:}{}{}
\theoremstyle{case}

\newcommand \rank{\mathrm{rank}}
\newcommand\MS{\operatorname{MinSupp}}
\newcommand \cl{\mathrm{cl}}

\usepackage[top=1.25in, bottom=1.25in, left=1in, right=1in]{geometry}                		         		

\usepackage{enumitem}
\usepackage{graphicx}
\usepackage{float}

\usetikzlibrary{decorations.pathmorphing}         

\newcommand{\com}[1]{\ignorespaces}

\begin{document}

\title{Matroids over skew tracts}

\author{Ting Su}

\email{tsu2@binghamton.edu}

\address{Department of Mathematics, Universität Hamburg, Germany}
\address{Present affiliation: Changjiang Geophysical Exploration and Testing Co., China}

\keywords{hyperfields, tracts, matroids, cryptomorphism}

\begin{abstract} Matroids over tracts provide an algebraic framework simultaneously generalizing the notions of matroids, oriented matroids, and valuated matroids, presented by Baker and Bowler. Pendavingh partially extended this theory to skew hyperfields and presented a new axiom system in terms of quasi-Pl\"ucker coordinates. We present a theory of matroids over skew tracts, which generalizes both the theory of matroids over tracts and the theory of weak matroids over skew hyperfields developed by Pendavingh. We give several cryptomorphic axiom systems for such matroids in terms of circuits, quasi-Pl\"ucker coordinates and dual pairs.
\end{abstract}

\maketitle

\section{Introduction}
%A \emph{matroid} (\cite{CR70}, \cite{Wel76}) is a combinatorial object that abstracts the notion of linear independence in a vector configuration over an arbitrary field. An \emph{oriented matroid} (\cite{BLV78}, \cite{LV75}) is a combinatorial object that abstracts the notion of linear independence in a vector configuration over an ordered field.The theories of matroids and oriented matroids are major branches of combinatorics with applications in many fields of mathematics, including topology, algebra, graph theory, and geometry. 

A \emph{matroid} can be thought of as a combinatorial object that abstracts the notion of linear independence in a vector configuration over an arbitrary field. An \emph{oriented matroid} can be thought of as a combinatorial object that abstracts the notion of linear independence in a vector configuration over an ordered field. 

In \cite{Bak16}, Baker and Bowler established {\em matroids over hyperfields}, which is an algebraic framework simultaneously generalizing the notion of linear subspaces, matroids, oriented matroids, valuated matroids, phased matroids, and some other ``matroids with extra structure". It gives a strong connection between combinatorial structures and algebraic structures; it also relates to tropical geometry.

A \emph{hyperfield} (\cite{Kra57}) is an algebraic structure similar to a field except that its addition is multivalued. Several good examples of hyperfields are given by Viro in \cite{Viro10}, including the \emph{Krasner hyperfield} $\mathbb{K}$, the \emph{sign hyperfield} $\mathbb{S}$ and the {\em ultratriangle hyperfield} $\mathbb{T}\triangle$. $\mathbb{T}\triangle$ is very important in tropical geometry. When $H$ is a field, an $H$-matroid corresponds to a linear subspace of some $H^n$. A $\mathbb{K}$-matroid is a usual matroid. An $\mathbb{S}$-matroid is an oriented matroid. And a $\mathbb{T}\triangle$-matroid is a valuated matroid, as defined in \cite{DW92b}. 

In \cite{Bak17}, Baker and Bowler defined a more general kind of algebraic object known as \emph{tracts}, which appear to be a natural setting for matroid theory. The other important example of tracts other than hyperfields is
given by \emph{partial fields} in the sense of Semple and Whittle (\cite{SW96}), which have also been the subject of much fruitful study.

Baker and Bowler also provided two natural notions of matroids over a tract $T$, {\em weak $T$-matroids} and {\em strong $T$-matroids}, which diverge for certain tracts. In particular, weak matroids and strong matroids are the same for fields, $\mathbb{K}$, $\mathbb{S}$ and $\mathbb{T}\triangle$. They provided both \emph{circuit axioms} and \emph{Grassman-Pl\"ucker axioms} for weak and strong $T$-matroids. The weak circuit axioms looks more like the \emph{(signed) circuit axioms} for (oriented) matroids and generalize them better than the strong axioms. The strong Grassmann-Pl\"ucker axioms looks more like the \emph{chirotope axioms} for oriented matroids and generalize them better than the weak axioms. Baker and Bowler also presented the dual pair axioms for both kinds of $T$-matroids and showed the cryptomorphism between the three axiom systems.

In \cite{Pen18}, Pendavingh partially extended the theory of matroids over tracts to {\em skew hyperfields}. He provided cryptomorphic axiom systems for {\em weak matroids over skew hyperfields} in terms of circuits and \emph{dual pairs}. However, there is no proper analogue of the Grassmann-Pl\"ucker function in the context of skew hyperfields. Pendavingh replaced the Grassmann-Pl\"ucker functions with \emph{quasi-Pl\"ucker coordinates} in the context of skew hyperfields. Then he provided a new axiom system for weak matroids over skew hyperfields in terms of quasi-Pl\"ucker coordinates. Moreover, Pendavingh also showed that all three axiom systems are cryptomorphic. 

In this paper, we introduce \emph{skew tracts}, which generalize tracts and skew hyperfields. We then present a theory of {\em matroids over skew tracts}, which generalizes the theory of matroids over tracts by Baker and Bower. Similarly, we provide two natural notions of matroids over a skew tract $T$, which we call {\em weak $T$-matroids} and {\em strong $T$-matroids}, corresponding to the two natural notions of matroids over tracts. 
The theory of weak matroids over skew tracts generalizes the theory of weak matroids over skew hyperfields by Pendavingh. But Pendavingh did not provide the theory of strong matroids over skew hyperfields and our theory of strong matroids over skew tracts fill this gap.
We also provide cryptomorphic axiom systems for both kinds of such matroids in terms of circuits, quasi-Pl\"ucker coordinates and dual pairs, and establish some basic duality results. 

\subsection{Structure of the paper}
We give the definitions of hyperfields and skew tracts in Section~\ref{sect.tract}. In Section~\ref{sect.matrodoverskewtract}, we provide different axioms for matroids over skew tracts in terms of circuits, dual pairs and quasi-Pl\"ucker coordinates and state the main results of duality theory. The rescaling and push-forward operations on matroids over skew tracts are also discussed in the same section. In Section~\ref{sect.crypt}, we present the proof for cryptomorphism of the three axiom systems.

\subsection{Acknowledgment} Thanks to Nathan Bowler for helpful discussion, comments and suggestions. Thanks also to Laura Anderson for comments and suggestions.

\section{Skew hyperfields and skew tracts}\label{sect.tract}

\begin{defn}\label{hyperoperation} A \textbf{hyperoperation} on a set $S$ is a map $\boxplus$ from $S \times S$ to the collection of non-empty subsets of $S$.

If $A$, $B$ are non-empty subsets of $S$, we define
$$A \boxplus B := \bigcup_{a\in A, b\in B} (a\boxplus b) $$
and we say that $\boxplus$ is \textbf{associative} if $a\boxplus(b\boxplus c)=(a\boxplus b) \boxplus c$ for all $a, b, c \in S$.
\end{defn}

All hyperoperations in this paper will be commutative and associative.

\begin{defn}\label{hypergroup} \cite{Bak16} 
A (commutative) \textbf{hypergroup} is a triple $(G,\boxplus, 0)$ where $\boxplus$ is a commutative and associative hyperoperation on $G$ such that:
\begin{enumerate}[(1)]
\item $0\boxplus x = \{x\}$ for all $x \in G$. 

\item For every $x \in G$ there is a unique element of G (denoted by $-x$ and called the {\bf hyperinverse} of $x$) such that $0 \in x \boxplus -x$.

\item (Reversibility) $x \in y \boxplus z$ if and only if $z \in x \boxplus -y$.
\end{enumerate}

A \textbf{skew hyperring} is a tuple $(R,\odot, \boxplus, 1, 0)$ such that:
\begin{enumerate}[(1)]
\item $(R, \odot, 1)$ is a monoid.

\item $(R, \boxplus, 0)$ is a hypergroup.

\item (Absorption rule) $x\odot 0 = 0\odot x = 0$ for all $x \in R$.

\item (Distributive Law) $a \odot (x \boxplus y )= (a \odot x) \boxplus (a \odot y)$ and $(x \boxplus y ) \odot a = (x \odot a) \boxplus (y \odot a)$ for all $a, x, y \in R$.
\end{enumerate}

A \textbf{hyperring} is a skew hyperring with commutative multiplication.

A skew hyperring $R$ is called a \textbf{skew hyperfield} if $0 \neq 1$ and every non-zero element of $R$ has a multiplicative inverse.

A \textbf{hyperfield} is then a skew hyperfield with commutative multiplication.
\end{defn}

%For $x, y \in F$, we will sometimes write $xy$ instead of $x\odot y$ if there is no risk of confusion.

%If $H$ is a skew hyperfield, then $H^E$ has a hypergroup structure given by 
%$$X \boxplus Y := \{Z \,|\, Z(e) \in X(e) \boxplus Y (e), \forall e\in E \}.$$
%$H$ acts on $H^E$ by componentwise multiplication.

\begin{example} Now we would like to introduce some more examples of skew hyperfields.
\begin{enumerate}
\item If $H$ is a skew field, then $H$ is a skew hyperfield with $x\odot y = x \cdot y$ and $x\boxplus y = \{x + y\}$, for $x, y\in H$.

\item The \textbf{Krasner hyperfield} $\mathbb{K}:=\{0, 1\}$ has the usual multiplication rule and hyperaddition is defined by $0\boxplus x = \{x\}$ for $x \in \mathbb{K}$ and $1\boxplus 1 = \{0,1\}$. 

\item The \textbf{sign hyperfield} $\mathbb{S}:= \{0, 1, -1\}$ has the usual multiplication rule and hyperaddition is defined by $0\boxplus x = \{x\}, x \boxplus x = \{x\}$ for $x \in \mathbb{S}$ and $1\boxplus -1 = \{0,1, -1\}$. 

\item The \textbf{ultratriangle hyperfield} $\mathbb{T}\triangle := \mathbb{R}_{\geq 0}$ (denoted by $\mathbb{Y}_{\times}$ in \cite{Viro10} and $\mathbb{T}$ in \cite{Bak17}) has the usual multiplication rule and hyperaddition is defined by
$$x \boxplus y =\begin{cases}
\{\max(x, y)\}	&\text{if $x \neq y$,}\\
\{z \,|\, z \leq x\} &\text{if $x = y$.}
\end{cases}
$$

\item The \textbf{phase hyperfield} $\mathbb{P}:= S^{1} \cup \{0\}$ has the usual multiplication rule and hyperaddition is defined by $0\boxplus x = \{x\}$ for $x \in \mathbb{P}$, $x \boxplus -x = \{0, x, -x\}$ and $ x \boxplus y = \{ \frac{a x + b y}{|a x + b y|} \,|\,  a, b \in \mathbb{R}_{>0} \}$ for $x, y \in S^{1}$ and $x\neq -y$. 

\item Let $G$ be a group with a binary operation $\cdot_G$. Let $H:= G\cup \{0\}$ with multiplication given by $x\cdot y = 0$ if $x$ or $y$ is 0 and by $\cdot_G$ of $G$ otherwise. Hyperaddition is given by
$$ x\boxplus y = \begin{cases}
H\backslash \{0\} & \text{ if } x \neq y. \\
H & \text{ if } x = y.
\end{cases}
$$
It is easy to check that $H$ is a skew hyperfield. When $G$ is abelian, $H$ is a hyperfield.

\item In \cite{Pen18}, Pendavingh showed that a field extension $K \subseteq L$ in positive characteristic $p$ and elements $x_e \in L$ for $e \in E$ gives rise to a matroid $M^\sigma$ on ground set $E$ with coefficients in a certain skew hyperfield $L^\sigma$. This skew hyperfield $L^\sigma$ is defined in terms of $L$ and its Frobenius action $\sigma: x \mapsto x^p$.
\end{enumerate}
\end{example}

\begin{defn}\label{tract} A \textbf{skew tract} is a group $G$ (written multiplicatively), together with an \textbf{additive relation structure} on $G$, which is a subset $N_G$ of the group semiring $\mathbb{N}[G]$ satisfying:
\begin{enumerate}[(1)]
\item The zero element of $\mathbb{N}[G]$ belongs to $N_G$.

\item The identity element $1$ of $G$ is not in $N_G$.

\item There is a unique element $\epsilon$ of $G$ with $1 + \epsilon \in N_G$.

\item $N_G$ is closed under the natural left and right actions of $G$ on $\mathbb{N}[G]$.
\end{enumerate}

A \textbf{tract} is a skew tract with the group $G$ abelian.
\end{defn}

One thinks of $N_G$ as those linear combinations of elements of G which ``sum to zero" (the $N$ in $N_G$ stands for ``null").

We let $T = G\cup \{0\}$ and $T^{\times} = G$. We often refer to the skew tract $(G, N_G)$ simply as $T$.

%It is tempting to associate a hyperadditive structure to $T$ by $a\boxplus b = \{ c\,|\, a+b-c \in N_G\}.$ However, when $T$ is a partial field, $a\boxplus b$ could be $\emptyset$ for some $a, b \in T$. We cannot define $\emptyset \boxplus a$, for $a\in T$. For instance, in $\mathbb{D}$, if we define $\emptyset \boxplus a = \emptyset$ for $a\in \mathbb{D}$, then$$(4\boxplus -1) \boxplus -1 = \emptyset \neq 4\boxplus (-1\boxplus -1) = 2.$$

Because of the following Lemma~\ref{-1}, we often write $-x$ instead of $\epsilon x$ for $x\in G$.

\begin{lem} \label{-1} Let $T$ be a skew tract.
\begin{enumerate}[(1)]
\item If $x, y \in G$ satisfy $x + y \in N_G$, then $y=\epsilon x = x \epsilon$.

\item $\epsilon^2 = 1$.

\item $G \cap N_G = \emptyset$.

\end{enumerate}
\end{lem}

For a skew hyperfield $H$ with hyperaddition $\boxplus$, we could define a skew tract $(G, N_G)$ associated to this skew hyperfield by setting $G = H\backslash \{0\}$ and $N_G = \{\sum_{i=1}^k x_i \,|\, \forall i, x_i \in G \text{ and } 0\in \boxplus_{i=1}^k x_i\}$.

If $T$ is a skew tract and $E$ is a non-empty finite set, we denote by $T^E$ the set of functions from $E$ to $T$. There are natural left and right actions of $T$ on $T^E$ by coordinate-wise multiplication. For $X, Y\in T^E$, we define
$$X+Y: = \{Z\in T^E\,|\, \forall e\in E, X(e)+Y(e) - Z(e) \in N_G\}.$$

The {\bf support} of $X\in T^E$, denoted by $\underline{X}$, is the set of $e\in E$ such that $X(e) \neq 0$.

\begin{defn}\label{involut}
Let $T$ be a skew tract. An \textbf{involution} of $T$ is a map $\tau: T \rightarrow T$, which preserves the addition, sends 0 to 0 and is a monoid homomorphism from $(T, \cdot, 1)$ to itself, such that $\tau^2$ is the identity map. We denote the image of $x \in T$ under $\tau$ by $\overline{x}$.
%An \textbf{involution} of $T$ is a homomorphism $\tau: T \rightarrow T$ such that $\tau^2$ is the identity map. Denote the image of $x \in T$ under $\tau$ by $\overline{x}$.
\end{defn}

\begin{defn}\label{orthg}
Let $T$ be a skew tract endowed with an involution $x \mapsto \overline{x}$. For $X, Y \in T^E$, the \textbf{product} of $X$ and $Y$ is defined as $$X\cdot Y := \sum_{e\in E} X(e)\cdot \overline{Y(e)}.$$
Note that $X \cdot Y \in \mathbb{N}[G]$. 
We say that $X, Y$ are \textbf{orthogonal}, denoted by $X \perp Y$, if $X \cdot Y \in N_G$. 

Let $\cC, \cD \subseteq T^E$. We say that $\cC, \cD$ are \textbf{orthogonal}, denoted by $\cC \perp \cD$, if $X \cdot Y \in N_G$ for all $X\in \cC$ and $Y\in \cD$. We say that $\cC, \cD$ are \textbf{$k$-orthogonal}, denoted by $\cC \perp_k \cD$, if $X \perp Y$ for all $X \in \cC$ and $Y \in \cD$ with $|\underline{X}\cap \underline{Y}|\leq k$.
\end{defn}

When $T$ is the field $\mathbb{C}$ of complex numbers or the phase matroid $\mathbb{P}$, the usual involution on $T$ is complex conjugation. For $T$ is $\mathbb{K}$, $\mathbb{S}$, or $\mathbb{T}\triangle$, the usual involution on $T$ is the identity map.

\section{Matroids over skew tracts} \label{sect.matrodoverskewtract}

In this section, we will define strong- and weak matroids over a skew tract $T$ or (for brevity) $T$-matroids. We will provide different axiom systems for both kinds of $T$-matroids in terms of circuits, dual pairs and quasi-Pl\"ucker coordinates.

\begin{nota} Throughout $E$ denotes a non-empty finite set. $T$ denotes a skew tract. For a skew tract $T$, $T^\times$ denotes $T-\{0\}$.

The $\bf{support}$ of $X \in T^E$ is $\underline{X} := \{e \in E \,|\, X(e) \neq 0\}.$ For $S \subseteq T^E$, $\supp(S)$ or $\underline{S}$ denotes the set of supports of elements of $S$, and $\MS(S)$ denotes the set of elements of $S$ of minimal support.

For simplicity, we usually write $E\backslash e$ for $E\backslash \{e\}$ when $e\in E$.
\end{nota}

We will always view a skew tract $T$ as being equipped with an involution $x \mapsto \overline{x}$. 

\subsection{Modular pairs}
First, we recall the definition of modular in the general context of lattices. %As a general reference on the combinatorics of posets and lattices we refer to \cite{Sta12}. 

%Let $(P, \leq)$ be a partially ordered set (poset). A \textbf{chain} in $P$ is a totally ordered subset $J$; the {\bf length} of a chain is $l(J):=|J|-1$. Given $x \in P$ we write $P_{\geq x} = \{x' \in P \,|\, x' \geq x\}$ and $P_{\leq x} = \{x' \in P \,|\, x' \leq x\}$. The {\bf length} of $P$ is $l(P):=\max \{l(J)\,|\, J \text{ a chain of } P \}$ and for $x \in P$, the {\bf height} of $x$ of $P$ is $l(x) := l(P_{\leq x})$.

%Let $x, y \in P$. If the poset $P_{\geq x}\cap P_{\geq y}$ has a unique minimal element, this element is denoted by $x\lor y$ and called the \textbf{join} of $x$ and $y$. If the poset $P_{\leq x}\cap P_{\leq y}$ has a unique maximal element, this element is denoted by $x\land y$ and called the \textbf{meet} of $x$ and $y$. The poset $P$ is called a \textbf{lattice} if the meet and join are defined for any $x, y \in P$. Every finite lattice $L$ has a unique minimal element $\hat{0}$ and a unique maximal element $\hat{1}$. 

\begin{defn}\cite{Oxl92, Bak17}\label{mod.latt} Let $L$ be a lattice with minimal element $\hat{0}$. An element $x\in L$ is called an \textbf{atom} if $x \neq \hat{0}$ and there is no $z \in L$ with $\hat{0} < z< x$. Two atoms $x, y \in L$ form a {\bf modular pair} if the height of $x \lor y$ is 2, i.e., $x\neq y$ and there do not exist $z_1,z_2\in L$ with $\hat{0} < z_1< z_2< x \lor y$. A family of atoms in $L$ is \textbf{modular} if the height of their join in $L$ is the same as the size of the family.
\end{defn}

Let $E$ be a set and let $C$ be a collection of pairwise incomparable nonempty subsets of $E$. The set $U(C):= \{\bigcup S\,|\, S \subseteq C \}$ forms a lattice when equipped with the partial order coming from inclusion of sets, with join corresponding to union and with the meet of $X$ and $Y$ defined to be the union of all sets in $C$ contained in both $X$ and $Y$. So every $X \in C$ is atomic as an element of $U(C)$. We say that $C_1, C_2 \in C$ form a \textbf{modular pair} in $C$ if they are a modular pair in $U(C)$; that is, the height of their join in the lattice $U(C)$ is 2. 

In a matroid $M$ of rank $r$ on $E$ with circuit set $C$, $C$ is a collection of pairwise incomparable nonempty subsets of $E$.

\begin{defn}\label{mod.MOM}
Two circuits $A, B \subseteq E$ form \textbf{a modular pair} in the matroid $M$ if $A\neq B$ and $A\cup B$ does not properly contain a union of two distinct elements of $C$. 
\end{defn}

Now we will define the modularity in $T$-matroids.

%\begin{defn} \cite{Del11, Bak17} Let $E$ be a set and let $C$ be a collection of pairwise incomparable nonempty subsets of $E$. We say that $C_1, C_2 \in C$ form a \textbf{modular pair} in $C$ if the height of their join in the lattice $U(C)$ is 2. 
%We say that $C_1, ... , C_k \in C$ form a \textbf{modular family} in $C$ if the height of their join in the lattice $U(C)$ is the same as the size of the family $k$. 
%\end{defn}

\begin{defn}\label{mod}
Let $\mathcal{C}$ be a subset of $T^E$. We say that $X, Y \in \cC$ form a \textbf{modular pair} in $\cC$ if $\underline{X}$, $\underline{Y}$ form a modular pair in $\supp(\cC)$. We say that $X_1, ... , X_k \in \cC$ form a \textbf{modular family} in $\cC$ if $\underline{X_1},... ,\underline{X_k}$ form a modular family in $\supp(\cC)$.
\end{defn}

\subsection{Circuit axioms} The following definition gives the circuit axioms for weak matroids over skew tracts. 

\begin{defn}\label{def.weakcir} Let $E$ be a non-empty finite set and let $T = (G, N_G)$ be a skew tract. A subset $\mathcal{C}$ of $T^E$ is called the \textbf{$T$-circuit set of a weak left $T$-matroid $\cM$ on $E$} if $\mathcal{C}$ satisfies the following axioms:
\begin{enumerate}[(C1)]
\item \label{c1} $0 \notin \mathcal{C}.$

\item \label{c2} (Symmetry) If $X \in \mathcal{C}$ and $\alpha \in T^{\times}$, then $\alpha \cdot X \in \mathcal{C}$.

\item \label{c3} (Incomparability) If $X$, $Y \in \mathcal{C}$ and $\underline{X} \subseteq \underline{Y}$, then there exists $\alpha \in T^{\times}$ such that $Y= \alpha \cdot X$.

\item \label{c4} (Modular Elimination) If $X$, $Y \in \mathcal{C}$ are a modular pair of $T$-circuits and $e\in E$ is such that $X(e)=-Y(e)\neq 0$, there exists a $T$-circuit $Z \in \mathcal{C}$ such that $Z(e)=0$ and $X(f) + Y(f) - Z(f)\in N_G $ for all $f \in E$.
\end{enumerate}
\end{defn}

We sometimes write $\cM$ as the ordered pair $(E, \cC)$. In the Modular Elimination axiom, we say $Z$ \textbf{eliminates} $e$ between $X$ and $Y$.

A {\bf weak right $T$-matroid} is defined analogously, with $\alpha \cdot X$ replaced by $X \cdot \alpha$ in (C2) and (C3). If $T$ is commutative, then weak left- and right $T$-matroids coincide, and we speak of weak $T$-matroids (refer to \cite{Bak17}).

If $\mathcal{C}$ is the set of $T$-circuits of a weak (left or right) $T$-matroid $\cM$ with ground set $E$, then there is an underlying matroid (in the usual sense) $\underline{\cM}$ on $E$ whose circuits are the supports of the $T$-circuits of $\cM$.

\begin{defn}\label{rank.T} The \textbf{rank} of $\cM$ is defined to be the rank of the underlying matroid $\underline{\cM}$.
\end{defn}

Now we present the circuit axioms for strong matroids over skew tracts.

\begin{defn}\label{def.strongcir}  A subset $\mathcal{C}$ of $T^E$ is called the \textbf{$T$-circuit set of a strong left $T$-matroid $\cM$ on $E$} if $\mathcal{C}$ satisfies \ref{c1}, \ref{c2} and \ref{c3} in Definition~\ref{def.weakcir} and the following stronger version of the Modular Elimination axiom \ref{c4}:
\begin{enumerate}[(C4)$'$]
\item \label{c4'} (Strong Modular Elimination) Suppose $X_1,... ,X_k$ and $X$ are $T$-circuits of $\mathcal{M}$ which together form a modular family of size $k+1$ such that $\underline{X} \not\subseteq \bigcup_{1\leq i \leq k}\underline{X_i}$, and for $1\leq i \leq k$ let
$$e_i \in (\underline{X}\cap \underline{X_i})\backslash \bigcup_{\substack{1\leq j \leq k \\ j\neq i} } \underline{X_j}$$
be such that $X(e_i)=-X_i(e_i)\neq 0$. Then there exists a $T$-circuit $Z \in \mathcal{C}$ such that $Z(e_i)=0$ for $1\leq i \leq k $ and $X(f)+ X_1(f)+ \cdots + X_k(f) - Z(f) \in N_G$ for every $f \in E$.
\end{enumerate}
\end{defn}

A {\bf strong right $T$-matroid} is defined analogously, with $\alpha \cdot X$ replaced by $X \cdot \alpha$ in (C2) and (C3). If $T$ is commutative, then strong left- and right $T$-matroids coincide, and we speak of strong $T$-matroids (refer to \cite{Bak17}).

From the definition, it is easy to see that any strong left (resp. right) $T$-matroid on $E$ is also a weak left (resp. right) $T$-matroid on $E$. 

A \textbf{projective $T$-circuit of a (weak or strong) left $T$-matroid} $\cM$ is an equivalence class of $T$-circuits of $\cM$ under the equivalence relation $X_1 \sim X_2$ if and only if $X_1 = \alpha \cdot X_2$ for some $\alpha \in T^\times$. Analogously, a \textbf{projective $T$-circuit of a (weak or strong) right $T$-matroid} $\cM$ is an equivalence class of $T$-circuits of $\cM$ under the equivalence relation $X_1 \sim X_2$ if and only if $X_1 = X_2 \cdot \alpha$ for some $\alpha \in T^\times$.

As \cite{Bak17} shows, a matroid over a field $T$ corresponds to a linear subspace of some $T^n$. A $\mathbb{K}$-matroid corresponds to a matroid. An $\mathbb{S}$-matroid is an oriented matroid. A $\mathbb{T}\triangle$-matroid is a valuated matroid. And a $\mathbb{P}$-matroid is a phased matroid, defined in \cite{AD12}. Weak matroids and strong matroids coincide over $\mathbb{K}$, $\mathbb{S}$ and $\mathbb{T}\triangle$, but they do not coincide over $\mathbb{P}$ (cf. Example 3.31 in \cite{Bak17}).

\subsection{Signature and coordinates} Now we will introduce two important terms for matroids.

\begin{defn} Let $N$ be a matroid on a finite set $E$ and let $T$ be a skew tract. We say a collection $\cC\subseteq T^E$ is a \textbf{left $T$-signature} of $N$ if $\cC$ satisfies \ref{c1}, \ref{c2} and \ref{c3} in Definition~\ref{def.weakcir}, and $\underline{\cC}$ is the set of circuits of the matroid $N$.
\end{defn}

Let $F$ be a subset of $E$ and let $a_1, ... , a_n$ be distinct elements of $E\backslash F$ with $n\in \mathbb{Z}_{>0}$. For simplicity, we write $Fa_1...a_n := F\cup \{a_1, ... , a_n\}$.

Let $N$ be a matroid of rank $r$ with bases $\cB$ and let $T$ be a skew tract. For $B, B'\in \cB$, we say $(B, B')$ is an ordered pair of {\bf adjacent bases} if $|B\backslash B'| = 1$. We name the set of ordered pairs of adjacent bases $A_N.$

\begin{defn}  A function $[\cdot] : A_N \rightarrow T$ comprises \textbf{left $T$-coordinates} for $N$ if $[\cdot]$ satisfies 
\begin{enumerate}[(LC1)]
\item \label{lc1} $[Fa, Fb]\cdot [Fb, Fa] = 1$\hfill if $|F| = r-1$ and $Fa, Fb \in \cB$.
\item \label{lc2} $[Fac, Fbc]\cdot [Fab, Fac]\cdot [Fbc, Fab] = -1$ \hfill if $|F| = r-2$ and $Fab, Fac, Fbc \in \cB$.
\item \label{lc3} $[Fac, Fbc] =  [Fad, Fbd]$ \hfill if $|F| = r-2$ and $Fac, Fad, Fbc, Fbd \in \cB$, but $Fab \notin \cB$.
\end{enumerate}
\end{defn}

Pendavingh showed in \cite{Pen18} that for a skew hyperfield $H$, a left $H$-signature encodes the same information as left $H$-coordinates. We will generalize this for skew tracts here.

Let $N$ be a matroid on $E$ and let $T$ be a skew tract. If $\cC$ is a left $T$-signature of $N$, then we may define a map $[\cdot]_{\cC}: A_N \rightarrow T$ by
$$[Fa, Fb]_{\cC} = -X(a)^{-1}X(b)$$
where $X\in \cC$ with $a, b\in \underline{X} \subseteq Fab$. This is well-defined, since if there exists $Y\in \cC$ such that $a, b\in \underline{Y} \subseteq Fab$, then $\underline{X} = \underline{Y}$ and thus by \ref{c3} there exists $\alpha\in T^\times$ so that $Y= \alpha \cdot X$. Then
$$Y(a)^{-1}Y(b) = (\alpha \cdot X(a))^{-1}\cdot (\alpha \cdot X(b)) = X(a)^{-1}X(b).$$
Conversely, given left coordinates $[\cdot]$ for $N$, we define
$$\cC_{N, [\cdot]} : = \{X\in T^E \,|\, \underline{X} \text{ is a circuit of } N \text{ and } X(a)^{-1}X(b) = -[Fa, Fb] \text{ whenever } a, b\in \underline{X} \subseteq Fab\}$$
We will usually omit the reference to $N$ when the choice of $N$ is unambiguous, and write $\cC_{[\cdot]}$.

\begin{lem}\label{sign.coord} Let $N$ be a matroid on $E$, let $T$ be a skew tract, let $[\cdot]: A_N \rightarrow T$, and let $\cC\subseteq T^E$. The following are equivalent.
\begin{enumerate}
\item $\cC$ is a left $T$-signature of $N$ and $[\cdot] = [\cdot]_{\cC}$.
\item $[\cdot]$ are left $T$-coordinates and $\cC = \cC_{[\cdot]}$.
\end{enumerate}
\end{lem}
The proof of the above lemma is the same as the proof of Lemma 1 of \cite{Pen18}.

The definition of {\bf right $T$-signatures} $\cC$, {\bf right $T$-coordinates} $[\cdot]$, and of the constructions $\cC_{[\cdot]}$ and $[\cdot]_{\cC}$ are obtained by reversing the order of multiplication throughout.

\subsection{Quasi-Pl\"ucker coordinates}
There is a cryptomorphic characterization of weak and strong matroids over a skew tract $T$ in terms of {\em quasi-Pl\"ucker coordinates}. We now describe it in this subsection.

\begin{defn} \label{def.weakQP} Let $T$ be a skew tract and let $N$ be a matroid of rank $r$ on $E$ with bases $\cB$. Then $[\cdot] : A_N \rightarrow T$ are \textbf{weak left quasi-Pl\"ucker coordinates} if
\begin{enumerate}[(P1)]
\item \label{p1} $[Fa, Fb]\cdot [Fb, Fa] = 1$\hfill if $|F| = r-1$ and $Fa, Fb \in \cB$.
\item \label{p2} $[Fac, Fbc]\cdot [Fab, Fac]\cdot [Fbc, Fab] = -1$ \hfill if $|F| = r-2$ and $Fab, Fac, Fbc \in \cB$.
\item \label{p3} $[Fa, Fb]\cdot [Fb, Fc]\cdot [Fc, Fa] = 1$ \hfill if $|F| = r-1$ and $Fa, Fb, Fc \in \cB$.
\item \label{p4} $[Fac, Fbc] =  [Fad, Fbd]$ \hfill if $|F| = r-2$, $Fac, Fad, Fbc, Fbd \in \cB$, and $Fab \notin \cB$ or $Fcd \notin \cB$.
\item \label{p5} $ -1 + [Fbd, Fab]\cdot [Fac, Fcd] + [Fad, Fab]\cdot [Fbc, Fcd] \in N_G$ 

\hfill if $|F| = r-2$ and $Fac$, $Fad$, $Fbc$, $Fbd$, $Fab$, $Fcd$ $\in \cB.$
\end{enumerate}
\end{defn}

It is clear to see from the definitions that weak left quasi-Pl\"ucker coordinates are left $T$-coordinates.

\begin{defn}\label{def.strongLQP} Let $r$ be the rank of $N$. $[\cdot] : A_N \rightarrow T$ are \textbf{strong left quasi-Pl\"ucker coordinates} if $[\cdot]$ satisfies \ref{p1}, \ref{p2} and \ref{p3} in Definition~\ref{def.weakQP} and the following stronger versions of \ref{p4} and \ref{p5} :
\begin{enumerate} [(P4)$'$]
\item \label{p4'} For any two subsets $I, J$ of $E$ with $|I| = r+1$, $|J| = r-1$ and $|I\backslash J|\geq 3$, we let $I_1 = \{x\in I \,|\, \text{both } I\backslash x \text{ and } Jx \text{ are bases of } N\}$. If $|I_1| = 2$ and we say $I_1 =\{a, b\}$, then $$[I\backslash a, I\backslash b] = [Jb, Ja].$$
\end{enumerate}
\begin{enumerate} [(P5)$'$]
\item \label{p5'} For any two subsets $I, J$ of $E$ with $|I| = r+1$, $|J| = r -1$ and $|I\backslash J|\geq 3$, we let $I_1 = \{x\in I \,|\, \text{both } I\backslash x \text{ and } Jx \text{ are bases of } N\}$. If $|I_1|\geq 3$, then for any $z\in I_1$,
$$-1 + \sum_{x\in I_1\backslash z} [I\backslash x, I\backslash z]\cdot [Jx, Jz] \in N_G.$$
\end{enumerate}
\end{defn}

It is easy to see that \ref{p4} is a special case of \ref{p4'} when $I = Fabd$ and $J = Fc$. \ref{p5} is also a special case of \ref{p5'} when $I = Fabd$, $J = Fc$ and $z = d$. 

The definition of {\bf weak right quasi-Pl\"ucker coordinates} and {\bf strong right quasi-Pl\"ucker coordinates} are obtained by reversing the order of multiplication throughout.

From the definition, any strong left (resp. right) quasi-Pl\"ucker coordinates are also weak left (resp. right) quasi-Pl\"ucker coordinates.

Quasi-P\"ucker coordinates are closely related to Grassmann Pl\"ucker functions for (commutative) hyperfields and Pendavingh showed this in Lemma 8 of \cite{Pen18}. The relationship will be the same for (commutative) tracts and we omit it in the paper.

\begin{thm} \label{crypt.qp.circ.weak}
Let $N$ be a matroid on $E$, let $T$ be a skew tract, let $[\cdot]: A_N \rightarrow T$, and let $\cC\subseteq T^E$. The following are equivalent.
\begin{enumerate}
\item $\cM = (E, \cC)$ is a weak left (resp. right) $T$-matroid such that $\underline{\cM} = N$ and $[\cdot] = [\cdot]_{\cC}$.
\item $[\cdot]$ are weak left (resp. right) quasi-Pl\"ucker coordinates for $N$ and $\cC = \cC_{[\cdot]}$.
\end{enumerate}
\end{thm}
The above theorem is proved by Theorem~\ref{crypt.all.weak}.

\begin{thm}\label{crypt.qp.circ.strong}
 Let $N$ be a matroid on $E$, let $T$ be a skew tract, let $[\cdot]: A_N \rightarrow T$, and let $\cC\subseteq T^E$. The following are equivalent.
\begin{enumerate}
\item $\cM = (E, \cC)$ is a strong left (resp. right) $T$-matroid such that $\underline{\cM} = N$ and $[\cdot] = [\cdot]_{\cC}$.
\item $[\cdot]$ are strong left (resp. right) quasi-Pl\"ucker coordinates for $N$ and $\cC = \cC_{[\cdot]}$.
\end{enumerate}
\end{thm}
The above theorem is proved by Theorem~\ref{crypt.all.strong}.

\subsection{Duality} In this subsection, we will present a duality theory for matroids over skew tracts which generalizes the established duality theory for matroids over tracts (cf. \cite{Bak17}).

\begin{thm}\label{dual} Let $E$ be a non-empty finite set with $|E|=n$, let $T$ be a skew tract, and let $\cM$ be a strong (resp. weak) left $T$-matroid of rank $r$ on $E$ with $T$-circuit set $\cC$ and strong (resp. weak) left quasi-Pl\"ucker coordinates $[\cdot]$. There is a strong (resp. weak) right $T$-matroid $\cM^*$ of rank $n-r$ on $E$, called the \textbf{dual matroid} of $\cM$, with the following properties:
\begin{enumerate}
\item Right $T$-coordinates $[\cdot]^*: A_{N^*}\rightarrow T$, defined by 
$$[B, B']^* : = -\overline{ [E\backslash B, E\backslash B']}$$
for all $(B, B') \in  A_{N^*}$, are strong (resp. weak) right quasi-Pl\"ucker coordinates of $\underline{\cM^*}$.
\item The $T$-circuit set of $\cM^*$ is $\cC^* := \MS(\{Y\in T^E\,|\, \forall X\in \cC, X\cdot Y \in N_G\}-\{\bf{0}\})$.
\item $\underline{\cM^*} = \underline{\cM}^*$.
\item $\cM^{**} = \cM$.
\end{enumerate}
If $\cM$ is a strong (resp. weak) right $T$-matroid, then the dual matroid $\cM^*$ is a strong (resp. weak) left $T$-matroid, $[\cdot]^*$ are strong (resp. weak) left quasi-Pl\"ucker coordinates, and the properties are obtained by reversing the order of multiplication throughout.

The $T$-circuits of $\cM^*$ are called the \textbf{$T$-cocircuits} of $\cM$, and vice-versa.
\end{thm}

\begin{proof}(cf. Theorem B in \cite{AD12}, Theorem 3.24 in \cite{Bak17})
This follows from Lemma~\ref{dual.qp}, Theorem~\ref{crypt.all.strong}, \ref{crypt.qp.dualpair}, and Lemma~\ref{dual.proof}.
\end{proof}

\subsection{Dual pairs} Now we will describe {\em dual pair} and present the dual pair axioms for matroids over skew tracts.

\begin{defn}\label{def.dual.pair.strong} Let $T$ be a skew tract, let $N$ be a matroid on $E$ and let $\cC, \cD \subseteq T^E$. We say that $(\cC, \cD)$ is a \textbf{dual pair of $T$-signatures} of $N$ if
\begin{enumerate}[(DP1)]
\item $\cC$ is a left $T$-signature of $N$, 
\item $\cD$ is a right $T$-signature of the dual matroid $N^*$,
\item $\cC \perp \cD.$
\end{enumerate}
\end{defn}

\begin{thm}\label{crypt.dual.circ.strong} Let $T$ be a skew tract, let $N$ be a matroid on $E$, let $\cC$ be a left $T$-signature of $N$, and let $\cD$ be a right $T$-signature of $N^*$. Then $\cC$ and $\cD$ are the sets of $T$-circuits and $T$-cocircuits, respectively, of a strong left $T$-matroid with underlying matroid $N$ if and only if $(\cC, \cD)$ forms a dual pair of $T$-signatures of $N$.
\end{thm}
The above theorem is proved by Theorem~\ref{crypt.all.strong}.

\begin{defn}\label{def.dual.pair.weak} We say that $(\cC, \cD)$ is a \textbf{weak dual pair of $T$-signatures} of $N$ if $\cC$ and $\cD$ satisfy (DP1), (DP2), and the following weakening of (DP3): 
\begin{enumerate} [(DP3)$'$]
\item $\cC \perp_3 \cD.$
\end{enumerate}
\end{defn}

\begin{thm}\label{crypt.dual.circ.weak} Let $T$ be a skew tract, let $N$ be a matroid on $E$, let $\cC$ be a left $T$-signature of $N$, and let $\cD$ be a right $T$-signature of $N^*$. Then $\cC$ and $\cD$ are the sets of $T$-circuits and $T$-cocircuits, respectively, of a weak left $T$-matroid with underlying matroid $N$ if and only if $(\cC, \cD)$ forms a weak dual pair of $T$-signatures of $N$.
\end{thm}
The above theorem is proved by Theorem~\ref{crypt.all.weak}.

It can be easily seen that the definition of $T$-matroids in terms of circuit axioms and quasi-Pl\"ucker coordinates axioms do not depend on the choice of the involution $\tau$ of $T$, however, the one in term of dual pairs does.

\subsection{Minors} Let $T$ be a skew tract, let $X\in T^E$, and let $A \subseteq E$. We define $X\backslash A\in T^{E\backslash A}$ by $(X\backslash A)(e)=X(e)$ for $e\in E\backslash A$. 

For $\cU\subseteq T^E$, define the {\bf deletion} of $A$ from $\cU$ as 
$$\cU\backslash A = \{ X\backslash A \,|\, X\in \cU, \underline{X}\cap A = \emptyset\}.$$
Define the {\bf contraction} of $A$ in $\cU$ as 
$$\cU/A = \MS(\{X\backslash A \,|\, X\in \cU\}).$$

\begin{thm}\label{minor} $\mathcal{C}\backslash A$ is the set of $T$-circuits of a strong (resp. weak) left $T$-matroid $\cM\backslash A$ on $E\backslash A$, called the \textbf{deletion} of $\cM$ by $A$, whose underlying matroid is $\underline{\cM}\backslash A$. 

Similarly, $\cC/A$ is the set of $T$-circuits of a strong (resp. weak) left $T$-matroid $\cM/A$ on $E\backslash A$, called the \textbf{contraction} of $\cM$ by $A$, whose underlying matroid is $\underline{\cM}/A$. 

Moreover, $(\cM\backslash A)^{*} = \cM^{*}/A$ and $(\cM/A)^{*} =\cM^{*}\backslash A $.
\end{thm}
The same statements with all instances of the word `left' replaced by `right' also hold.

\begin{proof} (cf. Theorem D in \cite{AD12}, Theorem 3.29 in \cite{Bak17}) This follows from Lemma~\ref{qp.del.cont}, Theorem~\ref{crypt.qp.dualpair} and \ref{crypt.all.strong}.
\end{proof}

\subsection{Strong and Weak matroids coincide over perfect skew tracts}
Given a (strong or weak) left matroid $\cM$ over a skew tract $T$ on a nonempty finite set $E$ with $T$-circuit set $\cC$ and $T$-cocircuit set $\cC^*$, a \textbf{$T$-vector} $V$ of $\cM$ is an element of $T^E$ such that for any $Y\in \cC^*$,
$V \cdot Y\in N_G.$
Similarly, a \textbf{$T$-covector} $U$ of $\cM$ is an element of $T^E$ such that for any $X\in \cC$,
$X\cdot U \in N_G.$ 
For a (strong or weak) right matroid, the $T$-vectors and $T$-covectors are obtained by reversing the order of multiplication throughout. We write $\cV(\cM), \cV^*(\cM)$ for the sets of $T$-vectors and $T$-covectors of $\cM$ respectively.

We say that a (strong or weak) left matroid $\cM$ is {\bf perfect} if $\cV(\cM)\perp \cV^*(\cM)$. A (strong or weak) right matroid $\cN$ is said to be {\bf perfect} if $\cV^*(\cN)\perp \cV(\cN)$.
We say that a skew tract $T$ is {\bf perfect} if for each strong left- and right $T$-matroid is perfect.

Skew fields, $\mathbb{K}$, $\mathbb{S}$, and $\mathbb{T}\triangle$ are perfect. But $\mathbb{P}$ is not perfect.

\begin{thm}
Any weak matroid $\cM$ over a perfect skew tract $T$ is strong.
\end{thm}
The proof of the above theorem is the same as the proof of Theorem 3.46 in \cite{Bak17}.

\subsection{Weak quasi-Pl\"ucker coordinates, duality and minors}

\begin{defn} Let $N$ be a matroid on $E$, let $T$ be a skew tract, and let $[\cdot]$ be weak left (resp. right) quasi-Pl\"ucker coordinates of $N$. Define the \textbf{dual map} $[\cdot]^*: A_{N^*} \rightarrow T$ of $[\cdot]$ by 
$$[B, B']^* : = - \overline{[E\backslash B, E\backslash B']}$$
for all $(B, B') \in  A_{N^*}$. It is evident that $[\cdot]^{**} = [\cdot]$.
\end{defn}

\begin{lem}\label{dual.qp} $[\cdot]^*$ are weak right (resp. left) quasi-Pl\"ucker coordinates of $N^*$. If $[\cdot]$ are strong left (resp. right) quasi-Pl\"ucker coordinates of $N$, then $[\cdot]^*$ are strong right (resp. left) quasi-Pl\"ucker coordinates of $N^*$.
\end{lem}
\begin{proof} We show the strong case first. First, we assume that $[\cdot]$ are strong left quasi-Pl\"ucker coordinates. We know that $\cB(N^*) = \{B \,|\, E - B \in \cB(N)\}$. Thus we just need to show that $[\cdot]^*$ satisfies \ref{p1}, \ref{p2}, \ref{p3}, \ref{p4'} and \ref{p5'} with order of multiplication reversed throughout.

\ref{p1}: Let $Fa, Fb$ be adjacent bases of $N^*$. Let $K= E\backslash (Fab)$. Thus $Ka, Kb$ are adjacent bases of $N$. Therefore by \ref{p1} of $[\cdot]$, we have
$$[Fb, Fa]^*\cdot [Fa, Fb]^* =  \overline{[Ka, Kb]}\cdot \overline{[Kb, Ka]} = 1.$$
So \ref{p1} holds for $[\cdot]^*$.

\ref{p2}: Let $Fab, Fbc, Fac$ be adjacent bases of $N^*$. Let $K= E\backslash (Fabc)$. Thus $Ka, Kb, Kc$ are adjacent bases of $N$. Therefore by \ref{p3} of $[\cdot]$, we have
\begin{align*}
&\,  [Fbc, Fab]^*\cdot [Fab, Fac]^* \cdot [Fac, Fbc]^*\\
= & \, - \overline{[Ka, Kc]} \cdot \overline{[Kc, Kb]}\cdot \overline{[Kb, Ka]}\\
= & \, -1.
\end{align*}
So \ref{p2} holds for $[\cdot]^*$.

\ref{p3}: Let $Fa, Fb, Fc$ be adjacent bases of $N^*$. Let $K= E\backslash (Fabc)$. Then $Kab, Kbc, Kac$ are adjacent bases of $N$. Therefore by \ref{p2} of $[\cdot]$, we have
$$[Fc, Fa]^* \cdot [Fb, Fc]^* \cdot [Fa, Fb]^* = - \overline{[Kab, Kbc]}  \cdot \overline{[Kac, Kab]}\cdot \overline{[Kbc, Kac]}=1.$$
So \ref{p3} holds for $[\cdot]^*$.

\ref{p4'} and \ref{p5'}: Let $r$ be the rank of $N$. Let $I, J$ be two subsets of $E$ with $|I| = r+1$, $|J|= r-1$, $|I\backslash J|\geq 3$, and let $I_1 = \{x\in I \,|\, \text{both } I\backslash x \text{ and } Jx \text{ are bases of } N^*\}$. Let $K = E\backslash I$ and let $L = E\backslash J$. Then $I\backslash x$ is a basis of $N^*$ if and only if $Kx$ is a basis of $N$. $Jx$ is a basis of $N^*$ if and only if $L\backslash x$ is a basis of $N$. So $I_1$ consists of every element $x$ such that both $Jx$ and $L\backslash x$ are bases of $N$. 

If $|I_1| = 2$, then we may say $I_1 = \{a, b\}$. Then by \ref{p4'} of $[\cdot]$ and Lemma~\ref{-1}, we have
$$[I\backslash a, I \backslash b]^* - [Jb, Ja]^* = - \overline{[Ka, Kb]} + \overline{[L\backslash b, L\backslash a]} \in N_G.$$
So $[I\backslash a, I \backslash b]^* = [Jb, Ja]^*$ and therefore \ref{p4'} holds for $[\cdot]^*$.

If $|I_1| \geq 3$, then we let $z\in I_1$. By \ref{p5'} of $[\cdot]$ we have
$$-1 + \sum_{x\in I_1\backslash z}  [Jx, Jz]^* \cdot [I\backslash x, I\backslash z]^*= -1+  \sum_{x\in I_1\backslash z} \overline{[L\backslash x, L\backslash z]} \cdot \overline{[Kx, Kz]} \in N_G.$$
Thus \ref{p5'} holds for $[\cdot]^*$. 

So by above argument, $[\cdot]^*$ are strong right quasi-Pl\"ucker coordinates of $N^*$.

Similarly, if $[\cdot]$ are strong right quasi-Pl\"ucker coordinates of $N$, then $[\cdot]^*$ are strong left quasi-Pl\"ucker coordinates of $N^*$.

The proof for weak case is essentially the same, but we only need to show \ref{p4} and \ref{p5} instead of \ref{p4'} and \ref{p5'}.
\end{proof}

Let $N$ be a matroid on $E$, let $T$ be a skew tract, let $[\cdot]: A_N \rightarrow T$ be weak (left or right) quasi-Pl\"ucker coordinates. Let $A$ be a subset of $E$, let $I_A\subseteq A$ be such that $I_A$ is a basis of $N|A$ and let $J_A\subseteq A$ be such that $J_A$ is a basis of $N/(E\backslash A)$. We know that $\cB(N/A) = \{B\,|\, B\cup I_A \in \cB(N)\}$ and $\cB(N\backslash A) = \{B \,|\, B\cup J_A\in \cB(N)\}$. So we can define the following two maps.

\begin{defn}\begin{enumerate}
\item (Contraction) Define $[\cdot]/A : A_{N/A} \rightarrow T$ by
$$[B, B']/A := [B\cup I_A, B'\cup I_A]$$
for all $(B, B') \in A_{N/A}.$
\item (Deletion) Define $[\cdot]\backslash A : A_{N\backslash A} \rightarrow T$ by
$$[B, B']\backslash A := [B\cup J_A, B'\cup J_A]$$
for all $(B, B') \in A_{N\backslash A}.$
\end{enumerate}
\end{defn}

\begin{lem}\label{qp.del.cont}\begin{enumerate}
\item $[\cdot]/A$ are weak left (resp. right) quasi-Pl\"ucker coordinates of $N/A$. If $[\cdot]$ are strong left (resp. right) quasi-Pl\"ucker coordinates, so is $[\cdot]/A$. The definition of $[\cdot]/A$ is independent of the choice of $I_A$.
\item $[\cdot]\backslash A$ are weak left (resp. right) quasi-Pl\"ucker coordinates of $N\backslash A$. If $[\cdot]$ are strong left (resp. right) quasi-Pl\"ucker coordinates, so is $[\cdot]\backslash A$. The definition of $[\cdot]\backslash A$ is independent of the choice of $J_A$.
\item $([\cdot]\backslash A)^* = [\cdot]^*/A.$
\end{enumerate}
\end{lem}
\begin{proof}
(1): The quasi-Pl\"ucker coordinates axioms are easy to check and we omit it here. Let $(B, B')\in A_{N/A}$ and let $I_A, I_A' \subseteq A$ be bases of $N|A$. To prove that the definition of $[\cdot]/A$ is independent of the choice of $I_A$, it suffices to prove that 
$$[B\cup I_A, B'\cup I_A] = [B\cup I_A', B'\cup I_A'].$$
We will prove the equation by induction on $|I_A \backslash I_A'|$. The base case $|I_A \backslash I_A'| = 1$ holds easily by \ref{p4}. Now let $k \geq 2$ and assume that the equation holds when $|I_A \backslash I_A'|\leq k-1$. So we may suppose that $|I_A \backslash I_A'| = k$.
Let $x\in I_A' \backslash I_A$. By basis exchange property, there exists $y \in I_A \backslash I_A'$ such that $(I_A' -x)\cup y =: I_A''$ is also a basis of $N|A$. Then we can see that $|I_A\backslash I_A''| = k-1$ and $|I_A''\backslash I_A'| = 1$. So by induction hypothesis and \ref{p4},
$$[B\cup I_A, B'\cup I_A] = [B\cup I_A'', B'\cup I_A''] = [B\cup I_A', B'\cup I_A'],$$
and we are done.

(2): The proof is similar as the proof in (1).

(3): similarly as the proof of Lemma 3.4 in \cite{AD12} we prove the case where $A = \{a\}$. Let $(B, B')\in A_{(N\backslash A)^*} = A_{N^*/A}$. Then let $I_A\subseteq A$ be such that $I_A$ is a basis of $N^*|A$ and let $J_A\subseteq A$ be such that $J_A$ is a basis of $N/(E\backslash A)$.

If $a$ is a coloop of $N$, then $a$ is a loop of $N^*$. So $I_A = \emptyset$, $J_A = \{a\}$, $a\notin B\cup B'$. Then
\begin{align*}
([B, B']\backslash a)^* & = - \overline{ [(E\backslash a)\backslash B, (E\backslash a)\backslash B']\backslash a }\\
& = - \overline{ [((E\backslash a)\backslash B)\cup J_A, ((E\backslash a)\backslash B')\cup J_A]} \\
& = - \overline{ [((E\backslash a)\backslash B)\cup a, ((E\backslash a)\backslash B')\cup a]} \\
& = - \overline{ [E\backslash B, E\backslash B']} \\
& = [B, B']^* \\
& = [B\cup I_A, B'\cup I_A]^* \\
& = ([B, B']^*)/a
\end{align*}

If $a$ is not a coloop of $N$, then $a$ is not a loop of $N^*$. So $I_A = \{a\}$, $J_A = \emptyset$. Then
\begin{align*}
([B, B']\backslash a)^* & = - \overline{ [(E\backslash a)\backslash B, (E\backslash a)\backslash B']\backslash a }\\
& = - \overline{ [((E\backslash a)\backslash B)\cup J_A, ((E\backslash a)\backslash B')\cup J_A] } \\
& = - \overline{ [(E\backslash a)\backslash B, (E\backslash a)\backslash B')] } \\
& = - \overline{ [E\backslash (B\cup a), E\backslash (B'\cup a)] } \\
& = [B\cup a, B' \cup a]^* \\
& = [B\cup I_A, B'\cup I_A]^* \\
& = ([B, B']^*)/a
\end{align*}
\end{proof}

\subsection{Homomorphism and push-forward} In this subsection, we will discuss homomorphisms of matroids over skew tracts and introduce an operation {\em push-forward}.

%\begin{defn}\label{hyp.homo}\cite{Bak16, Pen18} A \textbf{(commutative) hypergroup homomorphism} is a map $f : G \rightarrow H$ such that $f(0)= 0$ and $f(x \boxplus y) \subseteq f(x)\boxplus f(y)$ for all $x, y \in G$.

%A \textbf{skew hyperring homomorphism} is a map $f : R \rightarrow S$ which is a homomorphism of additive hypergroups as well as a homomorphism of multiplicative monoids (i.e., $f(1) = 1$ and $f(x\odot y) = f(x)\odot f(y)$ for $x, y \in R$).

%A \textbf{skew hyperfield homomorphism} is a homomorphism of the underlying skew hyperrings.
%\end{defn}

\begin{defn}\label{tract.homo} A \textbf{homomorphism} $f : (G, N_G) \rightarrow (H, N_H)$ of skew tracts is a group homomorphism $f : G \rightarrow H$, together with a map $f : \mathbb{N}[G] \rightarrow \mathbb{N}[H]$ satisfying $f(\sum_{i = 1}^k g_i) = \sum_{i=1}^k f(g_i)$ for $g_i \in G$, such that if $\sum_{i = 1}^k g_i \in N_G$ then $\sum_{i = 1}^k f(g_i) \in N_H$.
\end{defn}

An involution of a skew tract $T$ is a homomorphism from $T$ to itself.

\begin{lem} If $f : T \rightarrow T'$ is a homomorphism of skew tracts and $\cM$ is a strong (resp. weak) left $T$-matroid on set $E$ with $T$-circuit set $\cC$, then
$$\{\alpha \cdot f_* (X) \,|\, \alpha \in (T')^\times, X \in \cC\} =: f_*\cC$$
is the set of $T'$-circuit of a strong (resp. weak) left $T'$-matroid $f_*(\cM)$ on $E$, called the \textbf{push-forward} of $\cM$.
\end{lem}

From Definition~\ref{tract.homo}, it is not necessary that $f \circ \tau = \tau' \circ f$, where $\tau$ and $\tau'$ are involutions of $T$ and $T'$, respectively. So duality does not commute with the push-forward operation, that is, it is not necessary that $f_*(\cM^*) = (f_*(\cM))^*.$

It is easy to verify that the underlying matroids $\underline{f_*(\cM)} = \underline{\cM}$. Given Quasi-Pl\"ucker coordinates $[\cdot]_\cC: A_N \rightarrow T$ of $\underline{\cM}$, it is immediate that $[B, B']_{f_*\cC} = f([B, B']_\cC)$ for all adjacent bases $B$, $B'$. 

In particular, for every skew tract $T$, there is a unique homomorphism $\psi: T\rightarrow \mathbb{K}$ sending $0_T$ to $0_\mathbb{K}$ and sending every nonzero element of $T$ to $1_{\mathbb{K}}$. If $\cM$ is a $T$-matroid, then the push-forward $\psi_*(\cM)$ coincides with the underlying matroid $\underline{\cM}$.

\subsection{Rescaling} Now we will present an operation on matroids over skew tracts, called {\em rescaling}.

\begin{defn} Let $T$ be a skew tract, let $X\in T^E$ and let $\rho:E \rightarrow T^\times$. Then \textbf{right rescaling} $X$ by $\rho$ yields the vector $X\cdot\rho\in T^\times$ with entries $(X\cdot\rho)(e) = X(e)\cdot \rho(e)$ for all $e\in E$. Similarly, \textbf{left rescaling} gives a vector $\rho\cdot X$. We use $\rho^{-1}$ for the function from $E$ to $T^\times$ such that $\rho^{-1}(e) = \rho(e)^{-1}$ for all $e\in E$.

Let $\cC\subseteq T^E$. We define
$$\cC \cdot\rho:=\{X\cdot\rho \,|\, X\in \cC\} \text{ and } \rho\cdot\cC:=\{\rho\cdot Y\,|\, Y\in \cC\}.$$
\end{defn}

\begin{lem} Let $T$ be a skew tract, let $\cC$ be a left $T$-signature of a matroid $N$ on $E$, and let $\rho:E \rightarrow T^\times$. Then $\cC\cdot \rho$ is also a left $T$-signature of $N$.
\end{lem}
\begin{proof}
By definition of left $T$-signature, we know that $0\notin \cC$. Then $0\notin \cC\cdot \rho$. As $\underline{\cC\cdot \rho} = \underline{\cC}$, then $\underline{\cC\cdot \rho}$ is also the set of circuits of $N$. So it suffices to show that both (C2) and (C3) hold for $\cC\cdot \rho$.

(C2): Let $X\in \cC\cdot \rho$ and let $\alpha\in T^\times$. Then by definition, $X\cdot \rho^{-1} \in \cC$. By (C2) of $\cC$, we have that $(\alpha \cdot X) \cdot \rho^{-1} = \alpha \cdot (X\cdot \rho^{-1}) \in \cC$. So $\alpha \cdot X \in \cC\cdot \rho$.

(C3): Let $X, Y \in \cC\cdot \rho$ with $\underline{X} \subseteq \underline{Y}$. Then by definition $X\cdot \rho^{-1}, Y\cdot \rho^{-1} \in \cC$, and $\underline{X\cdot \rho^{-1}} \subseteq \underline{Y\cdot \rho^{-1}}$. By (C3) of $\cC$, there exists $\alpha\in T^\times$ such that $Y\cdot \rho^{-1} = \alpha \cdot (X\cdot \rho^{-1}) = (\alpha \cdot X) \cdot \rho^{-1}$. So $Y = \alpha \cdot X$.
\end{proof}

The same statement also holds with all instances of the word `left' replaced by `right' and $\cC\cdot \rho$ replaced by $\rho \cdot \cC$.

\begin{lem} Let $\cM$ be a strong (resp. weak) left $T$-matroid on $E$, let $\cN$ be a strong (resp. weak) right $T$-matroid on $E$, and let $\rho:E \rightarrow T^\times$.

Then $\cC(\cM)\cdot\rho^{-1}$ and $\rho\cdot\cC^*(\cM)$ are the $T$-circuit set and $T$-cocircuit set of a strong (resp. weak) left $T$-matroid $\cM^{\rho}$.
We say that $\cM^\rho$ arises from $\cM$ by \textbf{right rescaling}.

Similarly, $\rho^{-1}\cdot\cC(\cN)$ and $\cC^*(\cN)\cdot\rho$ are the $T$-circuit set and $T$-cocircuit set of a strong (resp. weak) right $T$-matroid $^\rho\cN$. We say that $^\rho\cN$ arises from $\cN$ by \textbf{left rescaling}.
\end{lem}

\begin{proof} By Theorem~\ref{crypt.dual.circ.strong} the statements hold for strong $T$-matroids because for any $\rho:E \rightarrow T^\times$, we have $\cC(\cM)\cdot\rho^{-1} \perp\rho\cdot\cC^*(\cM)$ if and only if $\cC(\cM) \perp \cC^*(\cM)$, and $\cC^*(\cN)\cdot\rho \perp \rho^{-1}\cdot\cC(\cN)$ if and only if $\cC^*(\cN) \perp \cC(\cN)$.

By Theorem~\ref{crypt.dual.circ.weak} the statements hold for weak $T$-matroids because for any $\rho:E \rightarrow T^\times$, we have $\cC(\cM)\cdot\rho^{-1} \perp_3\rho\cdot\cC^*(\cM)$ if and only if $\cC(\cM) \perp_3 \cC^*(\cM)$, and $\cC^*(\cN)\cdot\rho \perp_3 \rho^{-1}\cdot\cC(\cN)$ if and only if $\cC^*(\cN) \perp_3 \cC(\cN)$.
\end{proof}

\section{Cryptomorphism and duality} \label{sect.crypt}

\subsection{Between Quasi-Pl\"ucker coordinates and dual pairs}
This subsection will show that in the presence of an underlying matroid $N$, the quasi-Pl\"ucker coordinates axioms are cryptomorphic to the dual pair axioms.

\begin{lem}\label{2ortho.quasi} Let $N$ be a matroid on $E$, let $T$ be a skew tract, let $\cC \subseteq T^E$ be a left $T$-signature of $N$, and let $\cD\subseteq T^E$. The following are equivalent.
\begin{enumerate}
\item\label{2ortho.quasi.1} $\cD$ is a right $T$-signature of $N^*$ and $\cC \perp_2 \cD$.
\item $[\cdot] = [\cdot]_\cC$ satisfies \ref{p1}, \ref{p2}, \ref{p3} and \ref{p4'} and $\cD = \cC_{[\cdot]^*}$.
\end{enumerate}
\end{lem}
\begin{proof} The proof of this lemma with \ref{p4'} replaced by \ref{p4} is the same as the proof of Lemma 2 of \cite{Pen18}. As \ref{p4} is weaker than \ref{p4'}, we just need to show that  $[\cdot] = [\cdot]_\cC$ satisfies \ref{p4'} if (\ref{2ortho.quasi.1}) holds. So we assume that $\cD$ is a right $T$-signature of $N^*$ and $\cC \perp_2 \cD$. 

Let $r = \rank(N)$. Let $I, J$ be two subsets of $E$ with $|I| = r+1$, $|J| = r-1$ and $|I\backslash J|\geq 3$, and let $I_1 = \{x\in I \,|\, \text{both } I\backslash x \text{ and } Jx \text{ are bases of } N\}$ such that $|I_1|=2$. We suppose that $I_1 = \{a, b\}$. Then by definition, $I\backslash x$ is a basis of $N$ for $x\in I_1$. So there exists $X \in \cC$ such that $a, b \in \underline{X} \subseteq I$. By definition, $Jx$ is a basis of $N$ for all $x\in I_1$. So $E\backslash (Jx)$ is a basis of $N^*$ for $x\in I_1$. So there exists $Y \in \cD$ such that $a, b \in \underline{Y} \subseteq E\backslash J$. For any $z\in \underline{X}\backslash I_1$, $I\backslash z$ is a basis of $N$. So $Jz$ is dependent, and so $z\notin \underline{Y}$. So $\underline{X}\cap \underline{Y} = \{a, b\}$. By \ref{c2}, we may assume that $X(b) = 1 = Y(a)$. As $\cC\perp_2 \cD$, so $X\perp Y$. So 
$$N_G \ni X(a)\overline{Y(a)} + X(b)\overline{Y(b)} = X(a) + \overline{Y(b)}.$$
Then
$$[I\backslash a, I\backslash b] - [Jb, Ja] = [I\backslash a, I\backslash b] + \overline{[E\backslash (Jb), E\backslash (Ja)]^*} = -X(b)^{-1}X(a) - \overline{Y(b)Y(a)^{-1}} =- X(a) - \overline{Y(b)} \in N_G.$$

So $[I\backslash a, I\backslash b] = [Jb, Ja]$ and so \ref{p4'} holds. 
\end{proof}

\begin{thm}\label{crypt.qp.dualpair} Let $N$ be a matroid on $E$, let $T$ be a skew tract, let $\cC\subseteq T^E$ be a left $T$-signature of $N$, and let $\cD\subseteq T^E$. The following are equivalent.
\begin{enumerate}
\item $\cD$ is a right $T$-signature of $N^*$, and $(\cC, \cD)$ forms a strong (resp. weak) dual pair of $T$-signature of $N$, that is $\cC \perp \cD$ (resp. $\cC \perp_3 \cD$).
\item $[\cdot] = [\cdot]_\cC$ are strong (resp. weak) left quasi-Pl\"ucker coordinates and $\cD = \cC_{[\cdot]^*}$.
\end{enumerate}

Moreover, if (2) holds, then for any $e\in E$ we have
$$\cC_{[\cdot]\backslash e} = \cC \backslash e \text{ and }\cC_{[\cdot]/e} = \cC/e.$$
\end{thm}

\begin{proof} (cf. Lemma 3 in \cite{Pen18}) We will prove the strong case first. Let $n = |E|$ and $r= \rank(N)$. By Lemma~\ref{2ortho.quasi}, we only need to show that if $\cC$ is a left $T$-signature of $N$ and $\cD$ is a right $T$-signature of $N^*$ such that $\cC\perp_2 \cD$, then
$$\cC \perp \cD \text{ if and only if \ref{p5'} holds for } [\cdot] = [\cdot]_\cC.$$

We first show sufficiency. So we assume that $\cC \perp \cD$. Let $I, J$ be two subsets of $E$ with $|I| = r+1$, $|J| = r-1$ and $|I\backslash J|\geq 3$, and let $I_1 = \{x\in I \,|\, \text{both } I\backslash x \text{ and } Jx \text{ are bases of } N\}$ such that $|I_1|\geq 3$. We suppose that $I_1 = \{e_1, ... , e_k\}$ with $k\geq 3$. Then by definition, $I\backslash e_i$ is a basis of $N$ for all $1 \leq i\leq k$. So there exists $X \in \cC$ such that $I_1 \subseteq \underline{X} \subseteq I$. By definition, $Je_i$ is a basis of $N$ for all $1 \leq i\leq k$. So $E\backslash (Je_i)$ is a basis of $N^*$ for all $1 \leq i \leq k$. So there exists $Y \in \cD$ such that $I_1\subseteq \underline{Y} \subseteq E\backslash J$. For any $z\in \underline{X}\backslash I_1$, $I\backslash z$ is a basis of $N$. So $Jz$ is dependent, and so $z\notin \underline{Y}$. So $\underline{X}\cap \underline{Y} = I_1$. Let $z = e_1$. By \ref{c2}, we may assume that $X(e_1) = 1 = Y(e_1)$. As $\cC \perp \cD$, so $X\perp Y$. So
$$N_G \ni X(e_1)\overline{Y(e_1)} + ... + X(e_k)\overline{Y(e_k)} = 1 + \sum_{2 \leq i \leq k} X(e_i)\overline{Y(e_i)}.$$
Then
\begin{align*}
& -1 + \sum_{x \in I_1\backslash z} [I\backslash x, I \backslash z]\cdot [Jx, Jz] \\
= &  -1 - \sum_{x \in I_1\backslash z} [I\backslash x, I \backslash z]\cdot \overline{[E\backslash (Jx), E\backslash (Jz)]^*} \\
= & -1 - \sum_{x \in I_1\backslash z} X(z)^{-1}X(x)\cdot \overline{Y(x)Y(z)^{-1}} \\
= & -1 - \sum_{x \in I_1\backslash z} X(e_1)^{-1}X(x)\cdot \overline{Y(x)Y(e_1)^{-1}} \\
= & -1 - \sum_{x \in I_1\backslash z} X(x)\overline{Y(x)} \\
= &  -1 - \sum_{2\leq i \leq k} X(e_i)\overline{Y(e_i)} \\
\in & N_G.
\end{align*}
So \ref{p5'} holds.

We next show necessity. Let $X\in \cC$ and $Y\in \cD$ such that $|\underline{X}\cap \underline{Y}| \geq 3$. Let $K = \underline{X}\cap \underline{Y} = \{e_1, ... , e_k\}$. We extend $\underline{X}$ to a set $I$ of size $r+1$ such that $I\backslash e_i$ is a basis of $N$ for all $1 \leq i \leq k$. Similarly, we extend $\underline{Y}$ to a set $L$ of size $n-r+1$ such that $L\backslash e_i$ is a basis of $N^*$ for all $1 \leq i \leq k$. Let $J= E\backslash L$. We know that $I\backslash x$ is a basis of $N$ if and only if $x\in \underline{X}$. $Jy$ is a basis of $N$ if and only if $L\backslash y$ is a basis of $N^*$ if and only if $y\in \underline{Y}$. So $K$ consists of every element $x$ such that both $I\backslash x$ and $Jx$ are bases of $N$. By \ref{c2} we may assume that $X(e_1) = 1 = Y(e_1)$. By \ref{p5'} we have 
\begin{align*}
N_G & \ni -1 + \sum_{x\in K\backslash e_1} [I\backslash x, I\backslash e_1] \cdot [Jx, Je_1] \\
& = -1 - \sum_{x\in K\backslash e_1} [I\backslash x, I\backslash e_1] \cdot \overline{[L\backslash x, L\backslash e_1]^*} \\
& = - 1 - \sum_{x\in K\backslash e_1} X(e_1)^{-1}X(x)\cdot \overline{Y(x)Y(e_1)^{-1}}\\
& =  - X(e_1)\overline{Y(e_1)} - \sum_{x\in K\backslash e_1} X(x)\overline{Y(x)} \\
& = - \sum_{x\in K =\underline{X}\cap \underline{Y} }X(x)\overline{Y(x)} \\
& = - X\cdot Y.
\end{align*}
Thus $X\cdot Y \in N_G$ and therefore $X\perp Y.$

When (2) holds, we know that $\cC = \cC_{[\cdot]}$. So $\cC_{[\cdot]\backslash e} = \cC \backslash e$ and $\cC_{[\cdot]/e} = \cC/e$ follows immediately from the definition of $\cC$.

The proof for weak case is essentially the same, but in the special case that $|I_1|= k = 3$.
\end{proof}

So this shows the cryptomorphism between quasi-Pl\"ucker coordinates and dual pairs, and proves Theorem~\ref{dual} for dual matroids.

\begin{cor} \label{PivotProp} Let $N$ be a matroid on $E$, let $T$ be a skew tract, and let $[\cdot]: A_N \rightarrow T$ be weak left quasi-Pl\"ucker coordinates. We consider $X\in \cC_{[\cdot]}$ and $Y\in \cC_{[\cdot]^*}$. We extend $\underline{X}$ to $I$ such that $I\backslash x \in \cB(N)$ for all $x\in \underline{X}$. Similarly, we extend $\underline{Y}$ to $L$ such that $L\backslash y \in \cB(N^*)$ for all $y\in \underline{Y}$. Let $J = E\backslash L$. Then
\begin{enumerate}
\item (Pivoting property) for every $x_1, x_2 \in \underline{X}$, $$X(x_1)^{-1}X(x_2) = - [I \backslash x_2, I \backslash x_1] ;$$
\item  (Dual pivoting property) for every $y_1, y_2 \in \underline{Y}$, $$Y(y_1)Y(y_2)^{-1} = \overline{[ Jy_1, Jy_2]}.$$
\end{enumerate}
\end{cor}

% If $[\cdot]$ are right quasi-Pl\"ucker coordinates, then the properties will be given by reversing the order of multiplication on left hand side of the equations.

\subsection{From quasi-Pl\"ucker coordinates to circuits} In this subsection we will prove that the set $\cC_{[\cdot]}$ induced by quasi-Pl\"ucker coordinates $[\cdot]$ satisfies the circuit axioms.

\begin{thm}\label{crypt.qp.to.circ} Let $N$ be a matroid on $E$, let $T$ be a skew tract and let $[\cdot]: A_N \rightarrow T$ be strong (resp. weak) left quasi-Pl\"ucker coordinates for $N$. Then $\cM= (E, \cC_{[\cdot]})$ is a strong (resp. weak) left $T$-matroid such that $\underline{\cM} = N$ and $[\cdot] = [\cdot]_\cC$.
\end{thm}
\begin{proof} (cf. Theorem 4.13 in \cite{Bak17})
We prove the strong case first. Let $r = \rank(N)$. Let $[\cdot]$ be strong quasi-Pl\"ucker coordinates for $N$ and let $\cC = \cC_{[\cdot]}$. By Lemma~\ref{sign.coord}, we just need to show that $\cC$ satisfies \ref{c4'}.

Suppose we have a modular family $X, X_1, ... , X_k$ and elements $e_1, ... , e_k$ as in \ref{c4'} with $k\geq 2$. Let $z$ be any element of $\underline{X}\backslash \bigcup_{i=1}^k \underline{X_i}$. Let $A = \underline{X}\cup \bigcup_{i=1}^k \underline{X_i}$ and consider the matroid $N|A$. Since $A$ has height $k+1$ in the lattice of unions of circuits of $N$, then the rank of $N|A$ is $|A|- (k+1)$. Let $K = A\backslash \{z, e_1, ... , e_k\}$. The rank of $K$ is $|A|- (k+1)$, and so it is spanning in $N|A$. So $K$ is a basis of $N|A$. Let $O$ be a basis of $N/A$. Let $Z\in \cC$ with $\underline{Z}$ given by the fundamental circuit of $z$ with respect to $K\cupdot O$ and with $Z(z) = X(z)$. It is clear that $Z(e_i) = 0$ for all $1\leq i\leq k$. Now we would like to show that for any $f\in E$ we have
$$-Z(f) + X(f) + X_1(f) + \cdots + X_k(f) \in N_G.$$
If $f\notin A$ or $f\in \{z, e_1, ... , e_k\}$, then the statement holds directly. Thus we suppose that $f\in K$. Let $X_{k+1} \in \cC$ be such that $X_{k+1} = -Z$ and let $e_{k+1} = z$. So $X_{k+1}(e_{k+1}) = -X(e_{k+1})$ and $X_{k+1}(f) = -Z(f)$. Therefore we only need to show that for any $f\in K$ we have
$$X(f) + X_1(f) + \cdots + X_{k+1}(f) \in N_G.$$

We extend $\underline{X}\backslash e_{k+1}$ to $L$ such that $L$ is a basis of $N|A$. Then both $B_1 = L\cupdot O$ and $B_2 = K\cupdot O$ are bases of $N$. Let $I = \{e_{k+1}\}\cupdot B_1$ and let $J = B_2\backslash f$. Then $|I| = r+1$ and $|J| = r-1$. As $\{e_1, ..., e_{k+1}\}\subseteq \underline{X} \subseteq A$, then $\{e_1, ..., e_{k+1}\}\subseteq I\backslash J \subseteq \{e_1, ..., e_{k+1}, f\}$. Thus $|I\backslash J| \geq k+1 \geq 3.$

Let $I_1= \{x\in I \,|\, \text{both } I\backslash x \text{ and } Jx \text{ are bases of } N\} \subseteq I\backslash J$. So we only need to consider $x\in \{e_1, ..., e_{k+1}, f\}$. As $I\supseteq L\supseteq \underline{X}\backslash e_{k+1}$, then $I\backslash x$ is a basis of $N$ if and only if $x\in \underline{X}$.  For $1 \leq i \leq k+1$, we know that $K \supseteq \underline{X_i}\backslash e_i$, and so $e_i\in I_1$ if and only if $Je_i$ is a basis of $N$ if and only if $f\in \underline{X_i}$. As $Jf = B_2$ is a basis of $N$, then $f\in I_1$ if and only if $f\in \underline{X}$. As $f\in K\subseteq A$, then either $f\in \underline{X_i}$ for some $i\in \{1,...,k\}$ or $f\in \underline{X}$. Thus either $e_i \in I_1$ for some $i\in \{1,...,k\}$ or $f\in I_1$. By modularity of the family $X, X_1, ..., X_k$, $|I_1| \neq 1$.

\textbf{Case 1:} $|I_1| = 2$.

If $f\in \underline{X}$, then $I_1 = \{f, e_i\}$ for some $i\in \{1,...,k+1\}$. Moreover, $f\in \underline{X_i}$ and $f\notin \underline{X_j}$ for any $j\in \{1,...,k+1\}\backslash i$. By \ref{p4'}, we have $$[I\backslash f, I\backslash e_i] = [Je_i, Jf].$$
Thus $-X(e_i)^{-1}X(f) = -X_i(e_i)^{-1}X_i(f)$, and so $X(f) = -X_i(f)$. Therefore
$$X(f) + X_1(f) + \cdots + X_{k+1}(f) = X(f)+X_i(f) = -X_i(f)+X_i(f)\in N_G.$$

If $f\notin \underline{X}$, then $I_1 = \{e_i,e_j\}$ for some $i,j\in \{1,...,k+1\}$ with $i\neq j$. Moreover, $f\in \underline{X_i}$, $f\in \underline{X_j}$ and $f\notin \underline{X_l}$ for any $l\in \{1,...,k+1\}\backslash\{i,j\}$. By \ref{p4'}, \ref{p1} and \ref{p3}, we have $$[I\backslash e_i, I\backslash e_j] = [Je_j, Je_i] = [Je_j, Jf]\cdot [Jf, Je_i].$$
Thus $-X(e_j)^{-1}X(e_i) = X_j(e_j)^{-1}X_j(f)\cdot X_i(f)^{-1}X_i(e_i)$, and so $X_i(f)= -X_j(f)$. Therefore
$$X(f) + X_1(f) + \cdots + X_{k+1}(f) = X_i(f)+X_j(f) = -X_j(f) + X_j(f)\in N_G.$$

\textbf{Case 2:} $|I_1| \geq 3$.

If $f\in \underline{X}$, then $f\in I_1$ and $I_1\backslash f = \{e_i\,|\, f\in \underline{X_i} \text{ for } i\in \{1,...,k+1\}\}$. By \ref{p5'}, we have
\begin{align*}
N_G & \ni -1 + \sum_{x\in I_1\backslash f} [I\backslash x, I\backslash f]\cdot [Jx, Jf] \\
& = -1 + \sum_{\underset{f\in \underline{X_i}}{1\leq i\leq k+1}}  [I\backslash e_i, I\backslash f]\cdot [Je_i, Jf] \\
& = -1 + \sum_{\underset{f\in \underline{X_i}}{1\leq i\leq k+1}}  X(f)^{-1}X(e_i)\cdot X_i(e_i)^{-1}X_i(f) \\
& = -1 - \sum_{\underset{f\in \underline{X_i}}{1\leq i\leq k+1}}  X(f)^{-1}X_i(f) \\
& = -X(f)^{-1} \cdot (X(f) + \sum_{\underset{f\in \underline{X_i}}{1\leq i\leq k+1}}  X_i(f))\\
& = -X(f)^{-1} \cdot (X(f) + X_1(f) + \cdots +X_{k+1}(f)).
\end{align*}
So $$X(f) + X_1(f) + \cdots +X_{k+1}(f) \in N_G.$$

If $f\notin \underline{X}$, then $f\notin I_1$ and $e_l\in I_1$ for some $l\in \{1,...,k\}$. Without loss of generality, we assume that $l=1$. Then $I_1\backslash e_1 = \{e_i\,|\, f\in \underline{X_i} \text{ for } i\in \{2,...,k+1\}\}$. By \ref{p5'}, \ref{p1} and \ref{p3}, we have
\begin{align*}
N_G & \ni -1 + \sum_{x\in I_1\backslash e_1} [I\backslash x, I\backslash e_1]\cdot [Jx, Je_1] \\
& = -1 + \sum_{\underset{f\in \underline{X_i}}{2\leq i\leq k+1}}  [I\backslash e_i, I\backslash e_1]\cdot [Je_i, Je_1] \\
& = -1 + \sum_{\underset{f\in \underline{X_i}}{2\leq i\leq k+1}}  [I\backslash e_i, I\backslash e_1]\cdot [Je_i, Jf] \cdot [Jf, Je_1] \\
& = -1 + \sum_{\underset{f\in \underline{X_i}}{2\leq i\leq k+1}} -X(e_1)^{-1}X(e_i)\cdot X_i(e_i)^{-1}X_i(f)\cdot X_1(f)^{-1}X_1(e_1)\\
& =  -1 + \sum_{\underset{f\in \underline{X_i}}{2\leq i\leq k+1}} X(e_1)^{-1}X_i(f)\cdot X_1(f)^{-1}X_1(e_1)\\
& = X(e_1)^{-1}\cdot (X_1(f) + \sum_{\underset{f\in \underline{X_i}}{2\leq i\leq k+1}} X_i(f))\cdot X_1(f)^{-1}X_1(e_1)\\
& = X(e_1)^{-1}\cdot (X(f) + X_1(f) + \cdots + X_{k+1}(f))\cdot X_1(f)^{-1}X_1(e_1).
\end{align*}
So $$X(f) + X_1(f) + \cdots +X_{k+1}(f) \in N_G.$$

The proof for weak circuit axioms from weak quasi-Pl\"ucker coordinates axioms is essentially the same, but in the special case that $k=1$. This ensures that $|I\backslash J| \leq |\{z, f, e_1\}| = 3$, so that \ref{p4} and \ref{p5} can be applied instead of \ref{p4'} and \ref{p5'} respectively.
\end{proof}

\subsection{From circuits to dual pair} In this subsection, we will prove that the set $\cC$ of $T$-circuits of left (resp. right) $T$-matroid induces a (unique) right (resp. left) $T$-signature $\cD$ of the underlying matroid such that they form a dual pair of $T$-signature. We first introduce an important lemma.

\begin{lem} Let $T$ be a skew tract, let $\cM$ be a weak (left or right) $T$-matroid on a set $E$ and let $\cC$ be the set of $T$-circuits of $\cM$. Then for all $X, Y \in \cC$, $e, f\in E$ with $X(e) = -Y(e)\neq 0$ and $X(f)\neq - Y(f)$, there exists $Z\in \cC$ with $f\in \underline{Z}\subseteq \underline{X} \cup \underline{Y}\backslash e$.
\end{lem}
\begin{proof}
The proof of this lemma follows from the proof of Lemma 4.14 in \cite{Bak17}, the proof of Lemma 5.4 in \cite{AD12}.
\end{proof}

The following is a useful lemma from \cite{AD12}.

\begin{lem} \cite[Lemma 5.5]{AD12} Let $N$ be a matroid on the ground set $E$. Consider a circuit $C$ and a cocircuit $D$ of $N$ such that $|C \cap D| \geq 3$. Then there exist elements $e, f\in C \cap D$ and a cocircuit $D'$ of $N$ such that
\begin{enumerate}
\item $D$ and $D'$ are a modular pair,
\item $e\in C\cap D' \subseteq (C \cup D)\backslash f$.
\end{enumerate}
\end{lem}

\begin{thm}\label{crypt.circ.to.dualpair} Let $\cC$ be the $T$-circuit set of a weak left $T$-matroid $\cM$ on a set $E$. There is a unique right $T$-signature $\cD$ of $\underline{\cM}^*$ such that $(\cC, \cD)$ forms a weak dual pair of $T$-signatures of $\underline{\cM}$. If $\cM$ is a strong $T$-matroid, then $(\cC, \cD)$ forms a dual pair.
\end{thm}
\begin{proof} (cf. Proposition 5.6 in \cite{AD12}, Theorem 4.15 in \cite{Bak17}) We prove the strong case first. We assume that $\cM$ is a weak left $T$-matroid. For every cocircuit $D$ of $\underline{\cM}$, choose a maximal independent subset $A$ of $E\backslash D$. Then for every $e, f \in D$, there exists a unique circuit $C_{D,e,f}$ of $\underline{\cM}$ with support contained in $A\cup \{e,f\}$. By definition, there exists a unique projective $T$-circuit $X_{D,e,f}$ of $\cM$ with $\underline{X_{D,e,f}}= C_{D,e,f}$. For every circuit $X\in X_{D,e,f}$, it is easy to see that 
$$X(e)^{-1} X(f) = X_{D,e,f}(e)^{-1}X_{D,e,f}(f).$$

Now we define a collection $\cD \subseteq T^E$ by
$$\cD : = \{W \in T^E \,|\, D : = \underline{W} \in C^*(\underline{\cM}), \forall e, f\in \underline{W}, W(e)W(f)^{-1} = - \overline{X_{D,e,f}(e)^{-1}X_{D,e,f}(f)}\}.$$
We claim that $\cD$ is well defined. To show this claim, it is enough to prove that, given a cocircuit $D$ of $\underline{\cM}$ and $e, f, g\in D$,
$$-\overline{X_{D,e,f}(e)^{-1}X_{D,e,f}(f)} = ( -\overline{X_{D,e,g}(e)^{-1}X_{D,e,g}(g)}) \cdot (- \overline{X_{D,f,g}(g)^{-1}X_{D,f,g}(f)}).$$
The circuits $C_{D,e,g}$ and $C_{D,f,g}$ form a modular pair of $\underline{\cM}$, because their complements both contain the corank 2 coflat $\cl(E\backslash (A\cup \{e,f\}))$. Then (modular) elimination $g$ from $X_{D,e,g}(g)^{-1}X_{D,e,g}$ and $- X_{D,f,g}(g)^{-1}X_{D,f,g}$ gives $Y\in \cC$ with $e, f\in \underline{Y} \subseteq (\underline{X_{D,e,g}}\cup \underline{X_{D,f,g}})\backslash g \subseteq D\cup\{e,f\}$. So $Y\in X_{D,e,f}$. Moreover, $Y(e) = X_{D,e,g}(g)^{-1}X_{D,e,g}(e)$ and $Y(f) = - X_{D,f,g}(g)^{-1}X_{D,f,g}(f)$. So
\begin{align*}
& -X_{D,e,f}(e)^{-1}X_{D,e,f}(f) \\
= & -Y(e)^{-1}Y(f)\\
= &  - (X_{D,e,g}(g)^{-1}X_{D,e,g}(e))^{-1}\cdot (-X_{D,f,g}(g)^{-1}X_{D,f,g}(f)) \\
= &  (-X_{D,e,g}(e)^{-1}X_{D,e,g}(g)) \cdot (-X_{D,f,g}(g)^{-1}X_{D,f,g}(f))
\end{align*}
and so the claim follows.

Next we would like to show that $\cD$ is a right $T$-signature of $\cM$. By definition, we only need to prove that $\cD$ satisfies \ref{c2} and \ref{c3} with the order of multiplication reversed. 

For symmetry, let $W\in \cD$ and $\alpha\in T^\times$. Then for any $e, f\in \underline{W}=:D$, 
$$(W\cdot \alpha)(e)\cdot (W\cdot \alpha)(f)^{-1} = W(e)\cdot \alpha \cdot \alpha^{-1}\cdot W(f)^{-1} = W(e)W(f)^{-1} = - \overline{X_{D,e,f}(e)^{-1}X_{D,e,f}(f)}.$$
Thus $W\cdot \alpha \in \cD$.

For incomparability, let $W, Y\in \cD$ with $\underline{W} \subseteq \underline{Y}$. As $\underline{\cD} = C^*(\underline{\cM})$, then $\underline{W} = \underline{Y}$. Let $D:= \underline{W}\in C^*(\cM)$. Then for any $e, f\in D$,
$$W(e)W(f)^{-1} = - \overline{X_{D,e,f}(e)^{-1}X_{D,e,f}(f)} = Y(e)Y(f)^{-1}.$$
Thus $Y(e)^{-1}W(e) = Y(f)^{-1}W(f) = \alpha$ for some $\alpha\in T^\times$. Therefore $W(e) = Y(e)\cdot \alpha$ and $W(f) = Y(f)\cdot \alpha$. As $e, f$ are arbitrarily chosen, so $W = Y\cdot \alpha$.

Now it remains to prove that $\cC \perp \cD$. Let $X \in \cC$ and $Y \in \cD$. If $\underline{X} \cap \underline{Y}$ is empty, then we are done. So we suppose that $\underline{X} \cap \underline{Y}$ is nonempty. Since $\underline{\cM}$ is a matroid, $\underline{X} \cap \underline{Y}$ must contain at least two elements, so let $\underline{X} \cap \underline{Y} = \{z, e_1, ... , e_k\}$ with $k\geq 1$. By \ref{c2}, we may assume that $Y(z) = 1$. Let $K$ be a basis of $\underline{\cM}\backslash \underline{Y}$ including $\underline{X}\backslash \underline{Y}$. Then $B : = K\cup \{z\}$ is a basis of $\underline{\cM}$. For $1\leq i\leq k-1$, there exists $X_i\in \cC$ such that $\underline{X_i}$ is the fundamental circuit of $e_i$ with respect to $B$ in $\underline{\cM}$. By \ref{c2}, we may assume that $X_i(e_i) = -X(e_i)$ for all $1\leq i\leq k-1$. Let $C$ be the fundamental circuit of $e_k$ with respect to $B$ in $\underline{\cM}$.

We know that $\underline{X}\backslash B \subseteq \{e_1, ... , e_k\}$. For any $x\in \underline{X}\cap B$, the fundamental cocircuit of $x$ with respect to $E-B$ in $\underline{\cM}$ must meet $\underline{X}$ again, and must do so in some element of $\underline{X}\backslash B$. So $\underline{X} \subseteq C \cup \bigcup_{i=1}^{k-1} \underline{X_i}$, which has height $k$ in the lattice of unions of circuits of $\underline{\cM}$. So $X$ and the $X_i$ form a modular family of size $k$. So there exists $Z \in \cC$ with $Z(e_i) = 0$ for $1\leq i \leq k-1$ and $-Z(f) + X(f) + \sum_{i=1}^{k-1} X_i(f) \in N_G$ for all $f\in E$. Applying this with $f = e_k$ gives $Z(e_k) = X(e_k)$. By definition of $\cD$, we have $Y(e_i)Y(z)^{-1} = - \overline{X_i(e_i)^{-1}X_i(z)}$ for $1 \leq i \leq k-1$ and $Y(e_k)Y(z)^{-1} = -\overline{Z(e_k)^{-1}Z(z)}$. So for $1 \leq i \leq k-1$,
$$\overline{Y(e_i)} = -X_i(e_i)^{-1}X_i(z),$$
and
$$\overline{Y(e_k)} =-Z(e_k)^{-1}Z(z).$$
So \begin{align*}
X\cdot Y & = \sum_{x\in \underline{X}\cap \underline{Y}} X(x) \overline{Y(x)} \\
& = X(e_k) \overline{Y(e_k)} + X(z)\overline{Y(z)} + \sum_{i=1}^{k-1} X(e_i)\overline{Y(e_i)} \\
& =- X(e_k)\cdot Z(e_k)^{-1}Z(z) + X(z) - \sum_{i=1}^{k-1} X_(e_i)\cdot X_i(e_i)^{-1}X_i(z)\\
& =- Z(e_k)\cdot Z(e_k)^{-1}Z(z) + X(z) + \sum_{i=1}^{k-1} X_i(e_i)\cdot X_i(e_i)^{-1}X_i(z)\\
& = -Z(z) + X(z) + \sum_{i=1}^{k-1} X_i(z)\\
& \in N_G.
\end{align*}

The proof for weak dual pairs axioms from weak circuit axioms is essentially the same, but in the special case that $|\underline{X} \cap \underline{Y}| \leq 3$. %This ensures that $k = 1$ and $(X, X_1)$ forms a modular pair.
\end{proof}

The statement also holds when $\cM$ is a right $T$-matroid, $\cD$ will be defined by reversing the order of multiplication throughout, and $(\cD, \cC)$ is a dual pair.

\subsection{Cryptomorphic axiom systems for $T$-matroid} We can finally prove the main result of this paper. We begin by proving Theorems~\ref{crypt.qp.circ.strong} and \ref{crypt.dual.circ.strong} together in the following result:
\begin{thm}\label{crypt.all.strong}
Let $E$ be a finite set, let $N$ be a matroid and let $T$ be a skew tract. There are natural bijections between the following three kinds of objects:
\begin{enumerate}
\item Collections $\cC\subseteq T^E$ satisfying every axiom in Definition~\ref{def.strongcir}.
\item Maps $[\cdot]: A_N \rightarrow T$ satisfying every axiom in Definition~\ref{def.strongLQP}.
\item Dual pairs $(\cC, \cD)$ of $N$ satisfying every axiom in Definition~\ref{def.dual.pair.strong}.
\end{enumerate}
\end{thm}
\begin{proof} (1) $\Rightarrow$ (3): This is proved by Theorem~\ref{crypt.circ.to.dualpair}.

(3) $\Rightarrow$ (2): This is proved by strong case of Theorem~\ref{crypt.qp.dualpair}.

(2) $\Rightarrow$ (1): This is proved by Theorem~\ref{crypt.qp.to.circ}.
\end{proof}

Similarly we prove Theorems~\ref{crypt.qp.circ.weak} and \ref{crypt.dual.circ.weak} together in the following result:
\begin{thm}\label{crypt.all.weak}
Let $E$ be a finite set, let $N$ be a matroid and let $T$ be a skew tract. There are natural bijections between the following three kinds of objects:
\begin{enumerate}
\item Collections $\cC\subseteq T^E$ satisfying every axiom in Definition~\ref{def.weakcir}.
\item Maps $[\cdot]: A_N \rightarrow T$ satisfying every axiom in Definition~\ref{def.weakQP}.
\item A weak dual pair $(\cC, \cD)$ of $N$ satisfying every axiom in Definition~\ref{def.dual.pair.weak}
\end{enumerate}
\end{thm}

\subsection{Duality} Given the set $\cC$ of $T$-circuits of a (left or right) $T$-matroid, we will show that the corresponding set of $T$-cocircuits can be defined by orthogonality.

\begin{lem} \label{dual.proof}Let $T$ be a skew tract, let $\cM$ be a (strong or weak) left $T$-matroid on $T$, and let $\cC$ be the set of $T$-circuits of $\cM$. Then the set of elements of $\{Y\in T^E \,|\, X\perp Y, X\in \cC\} - \{0\}$ of minimal support is exactly the $T$-cocircuit set $\cD$ of $\cM$ given by Theorem~\ref{crypt.circ.to.dualpair}.
\end{lem}
\begin{proof}
The proof of this lemma is the same as the proof of Proposition 5.8 in \cite{AD12}.
\end{proof}

The statement also holds for (strong or weak) right $T$-matroids with  order of multiplication reversed throughout.

\bibliographystyle{alpha}
\bibliography{biblioCM}

\end{document}